\documentclass{article}
\usepackage{amssymb}
\newcommand\qed{\hbox{}\hfill\vrule width2mm height2mm}
\newcommand\F{{\cal F}}
\newcommand\hF{\hat{\cal F}}

\newcommand\M{{\cal M}}
\newcommand\A{{\cal A}}
\newtheorem{tet}{Theorem} 
\newtheorem{lem}{Lemma}
\newtheorem{cor}{Corollary}

\title{Largest family without $A\cup B\subseteq C\cap D$}
\author{{\bf Annalisa De Bonis} \\
University of Salerno
\\ Salerno, Italy \\E-mail: debonis@unisa.it 
\and
{\bf Gyula O.H. Katona}\thanks{The work of the second author was
supported by the Hungarian National Foundation for Scientific Research grant
number T037846, and UVO-ROSTE, Grant 875.630.9}
\\Alfr\'ed R\'enyi Institute of Mathematics, HAS,\\ Budapest P.O.B. 127 H-1364
HUNGARY \\ E-mail: ohkatona@renyi.hu
\and
{\bf Konrad J. Swanepoel }\thanks{The work of the third author was
supported by the South African National Research Foundation under Grant number 2053752.}
\\ Department of Mathematics, Applied Mathematics \\ and Astronomy   \\
University of South Africa \\ PO Box 392, UNISA 0003 
\\ South Africa \\
E-mail: swanekj@unisa.ac.za }

\begin{document} 
\date{}
\maketitle
\begin{abstract}
Let $\F$ be a family of subsets of an $n$-element set not containing four distinct members
such that $A\cup B\subseteq C\cap D$. It is proved that the maximum size of $\F$ under this 
condition is equal to the sum of the two largest binomial coefficients of order $n$. The 
maximum families are also characterized. A LYM-type inequality for such families is given, too.


\end{abstract}

\section{The inequalities}
Let $[n]=\{ 1,\ldots , n\} $ be a finite set and $\F\subset 2^{[n]}$ a family of its subsets.
The well-known theorem of {\it Sperner} (\cite{S}) says that if no member of $\F$ contains 
another member then $|\F |\leq {n\choose \lfloor {n\over 2}\rfloor}$, with equality iff $\F $ 
consists of all sets of size $\lfloor n/2\rfloor $ or all sets of size $\lceil n/2\rceil $. 
Moreover the LYM-type inequality (\cite{L} , \cite{Y} , \cite{M} , see also \cite{B})
$$\sum_{F\in \F}{n\choose |F|}^{-1}\leq 1$$ 
also holds for such a family. It is easy to see that the second inequality implies the first one.
On the other hand equality holds in the second inequality only when $\F$ consists of all sets
of a fixed size. 

The main aim of the present note is to investigate the analogous problem, when $\F$ contains
no four distinct sets $A, B, C, D$ such that $A$ is contained in both $C$ and $D$, and at the 
same time $B$ is contained in both $C$ and $D$. In other words, 
$${\mbox{there are no four distinct}}\ A,B,C,D\in \F \ {\mbox{ with}}\ A\cup 
B\subseteq C\cap D.\eqno(*)$$ 
It is easy to check that the family consisting of all $k$ and $k+1$-element subsets
satisfies ($*$). We will see that this is the largest family for the appropriate choice of $k$.
  
\begin{tet} Let $3\leq n$. If the family $\F \subseteq 2^{[n]}$ satisfies ($*$) then 
$|\F |$ cannot exceed the sum of the two largest binomial coefficients of order $n$, i.e.,
$|\F |\leq {n\choose \lfloor n/2\rfloor }+{n\choose  \lfloor n/2\rfloor +1}$.
\end{tet} 

The LYM-type inequality holds only if $\emptyset$ and $[n]$ are excluded from the family.

\begin{tet} Let $3\leq n$. If the family $\F \subseteq 2^{[n]}$ satisfies ($*$),
$\emptyset , [n]\not\in \F$ then
$$\sum_{F\in \F}{n\choose |F|}^{-1}\leq 2.$$
\end{tet}

Let us first prove Theorem 2 by the method of cyclic permutations (\cite{K}).
Let $\{ 1,\ldots ,n\} $ be considered as a cyclic permutation of the elements of $[n]$. 
That is the elements are considered $\pmod n$. An {\it interval}
is a subset of form  $\{ k,k+1,\ldots ,l\} $ where $1\leq k,l\leq n$. Intervals will be denoted by
$\hat{A}, \hat{B} $ etc. Families of intervals are denoted by $\hat{\cal A}, \hat{\cal B}$, etc.
The proof starts with two lemmas.

\begin{lem} Let $\hF$ be a family of intervals such that any member 
$\hat{F}\in \hF$ is contained in at most one other member of $\hF$, furthermore $\emptyset ,
[n]\not\in \hF$.
If $m$ denotes the number of the maximal members, $a$ denotes the number of non-maximal members 
then
$$m+{a\over 2}\leq n\eqno(1)$$ holds.
\end{lem}

{\bf Proof.} A {\it chain} is a family of subsets $L_1\subset \cdots \subset L_n$ where
$|L_i|=i (1\leq i\leq n)$. The number of chains of intervals containing $\hat{F}$ is 
$2^{|\hat{F}|-1}2^{n-|\hat{F}|-1}=2^{n-2}$. Suppose that 
$\hat{A}\subset \hat{B}, \hat{A}\not= \hat{B} $. We give an upper bound on the number of chains 
containing both of them. The number of choices of the new members of the chains 
``between the two sets'' is at most $2^{|\hat{B}|-|\hat{A}|-1}$ since, at least once, there is 
only one choice. Therefore the number of such chains is at most 
$2^{|\hat{A}|-1}2^{|\hat{B}|-|\hat{A}|-1}2^{n-|\hat{B}|-1}=2^{n-3}.$
The total number of chains is $n2^{n-2}$. Since a chain contains one or two members, we
obtain the inequality
$$(m+a)2^{n-2}\leq n2^{n-2}+a2^{n-3}$$
which is equivalent to the the statement of the lemma. \qed 
  
\begin{lem} 
If $\hat{\cal F}$ is a family of intervals satisfying ($*$), and $\emptyset , [n]\not\in \hat{\cal F}$, 
then $|\hat{\cal F}|\leq 2n$.
\end{lem}

{\bf Proof.} It is easy to see by complementation that the previous lemma holds for a family in 
which any member contains at most one other member. Divide $\hat{\cal F}$ into three subfamilies:
the maximal ($\hat{\cal M}_1$), the minimal ($\hat{\cal M}_2$) and other members ($\hat{\cal A}$). 
Introduce the notations $|\hat{\cal M}_1|=m_1, |\hat{\cal M}_2|=m_2, |\hat{\cal A}|=a$.
It is easy to see that ($*$) implies that $\hat{\cal M}_1\cup \hat{\cal A}$ satisfies the 
conditions 
of the previous lemma. Therefore we have $m_1+{a\over 2}\leq n$.
On the other hand $\hat{\cal M}_2\cup \hat{\cal A}$ 
satisfies the complementing of the previous 
lemma, we obtain the inequality $m_2+{a\over 2}\leq n$. The sum of the two inequalities is
$m_1+m_2+a\leq 2n$ as desired. \qed

{\bf Proof of Theorem 2.} We will double-count the pairs $({\cal C}, F)$ where ${\cal C}$ is 
a cyclic permutation of $[n]$, $F\in \F$ and $F$ is an interval along ${\cal C}$.  For a fixed 
$F$ the number of cyclic permutations is $|F|!(n-|F|)!$ therefore the number of pairs in 
question is 
$$\sum_{F\in \F}|F|!(n-|F|)! .$$
For a fixed cyclic permutation ${\cal C}$ the number of possible $F$s is at most $2n$ by the 
previous lemma. We obtained the the inequality
$$\sum_{F\in \F}|F|!(n-|F|)!\leq (n-1)!2n.$$
This is equivalent to the statement of the theorem. \qed

{\bf Proof of Theorem 1.} If none of $\emptyset$ and $[n]$ is a member of $\F$ then the 
statement is an easy consequence of Theorem 2. If both of them are in
$\F$ then $\F-\{ \emptyset , [n]\} $ is a Sperner family, therefore we have the upper estimate
${n\choose \lfloor n/2 \rfloor }+2$, which is less than our need, if $3\leq n$. Suppose that
exactly one of $\emptyset$ and $[n]$ is in $\F$. By complementation $\emptyset \in \F$ can be 
supposed. Then $\F ^{\prime}=\F -\{ \emptyset\} $ contains no 3 distinct members $A,B,C$ such that
$A\subset B, A\subset C$. It was proved in \cite{KT} (our Corollary 2 in Section 2 is slightly 
weaker) that 
$$|\F ^{\prime}|\leq \left( 1+{2\over n}\right) {n\choose \lfloor n/2 \rfloor }$$
holds under this condition. This upper estimate is strong enough when $3\leq n$. \qed

{\bf Remark.} D\'aniel Gerbner (student in Budapest) \cite{Ge}  noticed that there is no need to use
the theorem from \cite{KT}, since replacing $\emptyset$ by an arbitrarily chosen one-element set
$\{ i\} \not\in \F$ reduces the problem to the case when $\emptyset , [n]\not\in \F$. The case when
$\emptyset $ and all one-element sets are in $\F$ is trivial.
   
\section{Cases of equality}
The methods of the previous section are not strong enough for finding the cases of equality.
The conditions of Lemma 1 allow a large variety of families with equality. Therefore we have 
to consider the whole original family, rather than just the intervals. An {\it antichain} is a 
family of sets containg no comparable members.

\begin{lem} Let $\M$ and $\A$ be two disjoint antichains in $2^{[n]}$ where $[n]\not\in \M$. 
Suppose that for any $A\in \A$ there is a unique $f(A)\in \M$ with $A\subset f(A)$. Then
$$\sum_{M\in \M}{n\choose |M|}^{-1}+\sum_{A\in \A}{n\choose |A|}^{-1}
\left( 1-{1\over n-|A|}\right)\leq 1 \eqno(2)$$
holds, with equality only when either $|f(A)|=n-1$ or $|f(A)|=|A|-1$ holds for each $A\in \A$.
\end{lem}

{\bf Proof.} The number of chains containing a set $M$ is $|M|!(n-|M|)!$. Adding these  numbers
for all members of $\M$ and $\A$, a chain is counted once or twice, the latter can happen only
if the chain contains an $A\in \A$ and $f(A)\in M$. The total number of chains is $n!$, the number 
of chains containing both $A$ and $f(A)$ is $|A|!(|f(A)|-|A|)!(n-|f(A)|)!$.
Hence we have the following inequality:
$$\sum_{M\in \M}|M|!(n-|M|)!+  \sum_{A\in \A}|A|!(n-|A|)!\leq n! 
+ \sum_{A\in \A}|A|!(|f(A)|-|A|)!(n-|f(A)|)!.$$
Dividing by $n!$ we obtain
$$\sum_{M\in \M}{n\choose |M|}^{-1}+\sum_{A\in \A}{n\choose |A|}^{-1}
\left( 1-{n-|A|\choose n-|f(A)|}^{-1}\right)\leq 1. \eqno(3)$$
Since $|A|<|f(A)|<n$, the inequality $n-|A|\leq {n-|A|\choose n-|f(A)|}$ can be used in (3) 
to obtain (2). \qed

We know that $2\leq n-|A|$, which implies the following immediate corollary. 

\begin{cor} Under the conditions of Lemma 3 
$$\sum_{M\in \M}{n\choose |M|}^{-1}+\sum_{A\in \A}{1\over 2}{n\choose |A|}^{-1}
\leq 1 \eqno(4)$$
holds, with equality only when $|A|=n-2, |f(A)|=n-1$  for each $A\in \A$.
\end{cor}

\begin{cor} {\rm \cite{KT}} Let $4\leq n$. Suppose that the family $\F$ contains no three 
distinct members $A,B,C$ such that $A\subset B,C$. Then
$$|\F|\leq {n\choose \lfloor {n\over 2}\rfloor }\left( 1+{2\over n-3}\right)\eqno(5) $$
holds.
\end{cor}

{\bf Proof.} If $[n]\in \F$ then the rest of $\F$ satisfies the conditions of the Sperner theorem,
therefore we can suppose $[n]\not\in \F$. 
If we see that 
$${n\choose |A|}{n-|A|\over n-|A|-1}\leq {n\choose \lfloor {n\over 2}\rfloor }\left( 
1+{2\over n-3}\right)$$
holds for every $0\leq |A|\leq n-2$ then Lemma 3 implies (5). That is, we have to find the 
maximum of the function $g(i)={n\choose i}{n-i\over n-i-1}$ in the interval $0\leq i\leq n-2$.
Here $g(i-1)\leq g(i)$ holds if and only if $i(n-i-1)\leq (n-i)^2$. The discriminant 
$\sqrt{n^2-6n+1}$ of this quadratic inequality can be bounded from below and above by $n-4$ and
$n-3$, respectively, provided $7.5<n$. Hence $g(i-1)<g(i)$ holds if and only if 
$1\leq i< {n\over 2}+1$. The function $g(i)$ takes on its maximum 
in the interval $1\leq i\leq n-2$ at $\lfloor {n+1\over 2}\rfloor $.
The cases $n=4,5,6,7$ can be checked separately.   \qed

This corollary is slightly weaker than the statement in \cite{KT} , but its proof is much shorter.

\begin{tet} If $5\leq n$ then the equality in Theorems 1 and 2 implies that the 
family consist all $k$ and $k+1$-element sets for some $k$. For $n=3$ this is true only for Theorem 2,
in Theorem 1 there are other extremal constructions: the family
$$\{ \emptyset , \{ 1\} ,\{ 2\} , \{ 1,2 \} , \{ 2,3 \} , \{ 1,3\} \} $$
and its isomorphic versions. For $n=4$ there is up to isomorphism one more
extremal family for both theorems:
$$ {[4]\choose 2}\cup \{ \{ 1\} ,\{ 2,3,4\} ,\{ 2\} , \{ 1,3,4\} \} .$$  
\end{tet}

{\bf Proof.} Similarly to the proof of Lemma 2 define $\M_1$ and $\M_2$ as the families of 
maximal and minimal members of $\F$, respectively. $\A =\F -\M_1 -\M_2$. It is easy to see that
$\M_1 \cup \A$ satisfies the conditions of Corollary 1. On the other hand, the complements of
the members of $\M_2 \cup \A$ also satisfy it. The sum of the two inequalities again yield the statement
of Theorem 2. If $4<n$ there is no $A$ satisfying the conditions of equality in Corollary
1 for both (direct and complementing) cases. Therefore in this case the equality in Theorem
3 implies $\A =\emptyset$. It is well-known that $\F$ may consist of two full levels, only.
The cases $n=3,4$ can be checked separately. \qed

\end{document}